\documentclass[12pt]{article}

\title{\textbf{Epistemology as Information Theory:\\ From Leibniz to $\Omega$\footnote
{Alan Turing Lecture on Computing and Philosophy,
E-CAP'05,
European Computing and Philosophy Conference,
M\"alardalen University, 
V\"aster{\aa}s, Sweden, June 2005.}
}}

\author{\textbf{Gregory Chaitin}\thanks{IBM T. J. Watson Research Center, P. O. Box 218,
Yorktown Heights, NY 10598, U.S.A., \emph{chaitin@us.ibm.com.}}}

\date{}

\begin{document}

\maketitle

\begin{abstract}
In 1686 in his \emph{Discours de m\'etaphysique,} Leibniz points out
that if an arbitrarily complex theory is
permitted then the notion of ``theory'' becomes vacuous because there is
always a theory.  This idea is
developed in the modern theory of algorithmic information, which deals
with the size of computer programs
and provides a new view of G\"odel's work on incompleteness and Turing's
work on uncomputability.
Of particular interest is the halting probability $\Omega$, whose bits are
irreducible, i.e., maximally unknowable
mathematical facts. More generally, these ideas constitute a kind of
``digital philosophy'' related to recent attempts
of Edward Fredkin, Stephen Wolfram and others to view the world as a giant
computer. There are also connections
with recent ``digital physics'' speculations that the universe might
actually be discrete, not continuous. This
\emph{syst\`eme du monde} is presented as a coherent whole in my book \emph{Meta
Math!,} which will be published this fall.
\end{abstract}

\section*{Introduction}

I am happy to be here with you enjoying the delicate Scandinavian summer;
if we were a little farther north there wouldn't be any darkness at all.
And I am especially delighted to be here delivering the Alan Turing Lecture.
Turing's famous 1936 paper is an intellectual milestone that seems larger
and more important with every passing year.\footnote
{For Turing's original paper, with commentary, see Copeland's \emph{The Essential Turing.}}
    
People are not merely content to enjoy the beautiful summers in the far north,
they also want and need \textbf{to understand}, and so they create myths. In this part of
the world those myths involve Thor and Odin and the other Norse gods.
In this talk, I'm going to present another myth, what the French call a
\emph{syst\`eme du monde,} a system of the world, a speculative metaphysics based
on information and the computer.\footnote
{One reader's reaction (GDC): 
``Grand unified theories may be like myths, but surely there is a difference
between scientific theory and any other narrative?''
I would argue that a scientific narrative is more successful than the Norse myths
because it explains what it explains more precisely and without having to postulate
new gods all the time, i.e., it's a better ``compression'' 
(which will be my main point in this lecture;
that's how you measure how successful a theory is).}
    
The previous century had logical positivism and all that emphasis on the philosophy of language,
and completely shunned speculative metaphysics, but a number of us think that it is time
to start again. 
There is an emerging digital philosophy and digital physics, a new metaphysics associated with 
names like Edward Fredkin and Stephen Wolfram and a handful of like-minded individuals,
among whom I include myself. As far as I know the terms ``digital philosophy'' and
``digital physics'' were actually invented by Fredkin, and he has a large website
with his papers and a draft of a book about this.  Stephen Wolfram attracted a great
deal of attention to the movement and stirred up quite a bit of controversy with his very
large and idiosyncratic book on \emph{A New Kind of Science.}
    
And I have my own book on the subject, in which I've attempted to wrap everything I know
and care about into a single package. It's a small book, and amazingly enough it's
going to be published by a major New York publisher a few months from now.
This talk will be an overview of my book, which presents my own personal version of
``digital philosophy,'' since each of us who works in this area 
has a different vision of this tentative, emerging world view.
My book is called \emph{Meta Math!,} which may not seem like a serious title, but it's
actually a book intended for my professional colleagues as well as for the general public,
the high-level, intellectual, thinking public.
    
``Digital philosophy'' is actually a neo-Pythagorean vision of the world,
it's just a new version of that. According to Pythagoras, all is number ---
and by number he means the positive integers, 1, 2, 3, \ldots\
--- and God is a mathematician.
``Digital philosophy'' updates this as follows: 
Now everything is made out of 0/1 bits,
everything is digital software, and God is a computer programmer, not a mathematician!
It will be interesting to see how well this vision of the world succeeds,
and just how much of our experience and theorizing can be included or shoe-horned within this
new viewpoint.\footnote
{Of course, a system of the world can only work by omitting everything that doesn't
fit within its vision. The question is how much will fail to fit, and conversely,
how many things will this vision be able to help us to understand.  Remember, if one is wearing
rose colored glasses, everything seems pink.   And as Picasso said, theories are
lies that help us to see the truth. No theory is perfect, and it will be interesting
to see how far this digital vision of the world will be able to go.}
    
Let me return now to Turing's famous 1936 paper.
This paper is usually remembered for inventing the programmable digital
computer via a mathematical model, the Turing machine,
and for discovering the extremely fundamental halting problem.
Actually Turing's paper is called ``On computable numbers, with an application
to the \emph{Entscheidungsproblem,}'' and by computable numbers Turing means ``real''
numbers, numbers like $e$ or $\pi = 3.1415926\ldots$ that are measured with
infinite precision, and that can be computed with arbitrarily high precision, digit by digit 
without ever stopping,
on a computer.
    
Why do I think that Turing's paper ``On computable numbers'' is so important? Well, in my opinion
it's a paper on epistemology, because we only understand something if we
can program it, as I will explain in more detail later.  And it's a paper on physics,
because what we can actually compute depends on the laws of physics in 
our particular universe and distinguishes it from other possible universes.
And it's a paper on ontology, because it shows that some real numbers are
\textbf{uncomputable}, which I shall argue calls into question their very existence,
their mathematical and physical existence.\footnote
{You might exclaim (GDC), ``You can't be saying that before Turing and the computer
no one understood anything; that can't be right!''
My response to this is that
before Turing (and my theory) people could understand things,
but \textbf{they couldn't measure how well} they understood them.
Now you can measure that, in terms of the degree of compression that is achieved.
I will explain this later at the beginning of the section on computer epistemology.
Furthermore, programming something forces
you to understand it better, 
it forces you to really understand it,
since you are explaining it \textbf{to a machine}.
That's sort of what happens when a student or a small child asks you what at first
you take to be a stupid question, and then you realize that this question has in
fact done you the favor of forcing you to formulate your ideas more clearly and
perhaps even question some of your tacit assumptions.}
    
To show how strange uncomputable real numbers can be, let me give a particularly
illuminating example of one, which actually preceded Turing's 1936 paper.
It's a very strange number that was invented in a 1927 paper by the French mathematician Emile Borel. 
Borel's number is sort of an anticipation, a partial anticipation, of Turing's 1936
paper, but that's only something that one can realize in retrospect. 
Borel presages
Turing, which does not in any way lessen Turing's important contribution    
that so dramatically and sharply clarified all these vague ideas.\footnote
{I learnt of Borel's number by reading Tasic's \emph{Mathematics and the Roots of
Postmodern Thought,} which also deals with many of the issues discussed here.}
    
Borel was interested in ``constructive'' mathematics, 
in what you can actually compute we would say nowadays.
And he came up with an extremely strange non-constructive
real number.  You list all possible yes/no questions in French in an immense, an
infinite list of all possibilities.  This will be what mathematicians call a denumerable
or a countable infinity of questions, because it can be put into a one-to-one correspondence
with the list of positive integers 1, 2, 3, \ldots\  In other words, there will be a first
question, a second question, a third question, and in general an $N$th question.
    
You can imagine all the possible questions to be ordered by size, and within questions of
the same size, in alphabetical order.  More precisely, you consider
all possible strings, all possible finite sequences of symbols in the French alphabet,
including the blank so that you get words, and the period so that you have sentences.  
And you imagine filtering out
all the garbage and being left only with grammatical yes/no questions in French.
Later I will tell you in more detail how to actually do this. 
Anyway, for now \textbf{imagine} doing this,
and so there will be a first question, a second question, an $N$th question.
    
And the $N$th  digit or the $N$th bit after the decimal point of Borel's number answers
the $N$th question: It will be a 0 if the answer is no, and it'll be a 1 if the answer is yes.
So the binary expansion of Borel's number contains the answer to every possible yes/no 
question!
It's like having an oracle, a Delphic oracle that will answer every yes/no question!
    
How is this possible?!  Well, according to Borel, it isn't really possible, this can't be,
it's totally unbelievable.
This number is only a mathematical fantasy, it's not for real, it cannot claim a legitimate
place in our ontology.
Later I'll show you a modern version of Borel's number, my halting probability $\Omega$.
And I'll tell you why some contemporary physicists, real physicists, not mavericks, are
moving in the direction of digital physics.

{\small
[Actually, to make Borel's number as real as possible, you have to avoid the problem of filtering
out all the yes/no questions. And you have to use decimal digits, you can't use binary digits.
You number all the possible finite strings of French symbols including blanks and periods, which is
quite easy to do using a computer.
Then the $N$th digit of Borel's number is 0
if the $N$th string of characters in French is ungrammatical and not proper French,
it's 1 if it's grammatical, but not a yes/no question, 
it's 2 if it's a yes/no question that cannot be answered (e.g., ``Is the answer to
this question `no'?''), it's 3 if the answer is no, and it's 4 if the answer is yes.]
}
    
\emph{Geometrically} a real number is the most straightforward thing in the world, it's
just a point on a line. That's quite natural and intuitive. But \emph{arithmetically,}
that's another matter. 
The situation is quite different.
From an arithmetical point of view reals are extremely problematical, 
they are fraught with difficulties!
     
Before discussing my $\Omega$ number, 
I want to return to the fundamental question of what does it mean to understand.
How do we explain or comprehend something? What is a theory? 
How can we tell whether or not it's a successful theory?
How can we measure how
successful it is?
Well, using the ideas of information and computation, that's not difficult to do,
and the central idea can even be traced back 
to Leibniz's 1686 \emph{Discours de m\'etaphysique.}

\section*{Computer Epistemology: What is a mathematical or scientific theory?
How can we judge whether it works or not?}

In Sections V and VI of his \emph{Discourse on Metaphysics,}
Leibniz asserts that God simultaneously maximizes the variety, diversity and richness of the 
world, and minimizes the conceptual complexity of the set of ideas that determine the world.
And he points out that for any finite set of points there is always a mathematical equation
that goes through them, in other words, a law that determines their positions.
But if the points are chosen at random, that equation will be extremely complex.
     
This theme is taken up again in
1932 by Hermann Weyl in his book \emph{The Open World} consisting of 
three lectures he gave at Yale University on the
metaphysics of modern science.
Weyl formulates 
Leibniz's crucial idea in 
the following extremely dramatic fashion: If one permits
arbitrarily complex laws, then the concept of law becomes vacuous, because there is always a law! 
Then Weyl asks, how can we make more precise the distinction between 
mathematical simplicity and mathematical complexity? It seems to be very hard to do that.
How can we measure this important parameter, without which it is impossible
to distinguish between a successful theory and one that is completely unsuccessful?
     
This problem is taken up and I think satisfactorily resolved in the new mathematical theory
I call \emph{algorithmic information theory.} The epistemological model that is central to
this theory is that a scientific or mathematical theory 
is a computer program for calculating the facts, and the
smaller the program, the better.  The complexity of your theory, of your law, is measured
in bits of software:  
\begin{center}
program (bit string) $\longrightarrow$ \textbf{Computer} $\longrightarrow$ output (bit string)
\\ 
theory $\longrightarrow$ \textbf{Computer} $\longrightarrow$ mathematical or scientific facts
\\  
Understanding is compression!
\end{center}
     
Now Leibniz's crucial observation can be formulated much more precisely.  For any finite set
of scientific or mathematical facts, there is always a theory that is exactly as
complicated, exactly the same size in bits, as the facts themselves.  
(It just directly outputs them ``as is,'' without doing any computation.)
But that doesn't count, that doesn't enable us to distinguish between what can be comprehended
and what cannot, because there is always a theory that is as complicated as what it explains.
A theory, an explanation, is only successful to the extent to which it compresses the number of
bits in the facts into a much smaller number of bits of theory.  Understanding is compression,
comprehension is compression!
That's how we can tell the difference between real theories and \emph{ad hoc} theories.\footnote
{By the way, Leibniz also mentions complexity in 
Section 7 of
his \emph{Principles of Nature and Grace,}
where he asks the amazing question, ``Why is there something rather than nothing? For nothing
is simpler and easier than something.''}
    
What can we do with this idea that an explanation has to be simpler than what it explains?  
Well, the most important application of these ideas that I
have been able to find is in metamathematics, it's in discussing what mathematics can or cannot
achieve. You simultaneously get an information-theoretic, computational perspective on G\"odel's
famous 1931 incompleteness theorem, and on Turing's famous 1936 halting problem.  How?\footnote
{For an insightful treatment of G\"odel as a philosopher, see Rebecca Goldstein's
\emph{Incompleteness.}}
    
Here's how! These are my two favorite information-theoretic incompleteness results:
\begin{itemize}
\item
You need an $N$-bit theory in order to be able to prove that a specific $N$-bit program
is ``elegant.''
\item
You need an $N$-bit theory in order to be able to determine $N$ bits of the 
numerical value, of the base-two binary expansion, of the halting probability $\Omega$.
\end{itemize}
Let me explain.
    
What is an elegant program?  It's a program with the property that no program written in
the same programming language that produces the same output is smaller than it is. In other 
words, an elegant program is the most concise, the simplest, the best theory for its output.
And there are infinitely many such programs, they can be arbitrarily big, because for any
computational task there has to be at least one elegant program. (There may be several
if there are ties, if there are several programs for the same output that have exactly
the minimum possible number of bits.)
    
And what is the halting probability $\Omega$? Well, it's defined to be the probability
that a computer program generated at random, by choosing each of its bits using an
independent toss of a fair coin, will eventually halt.
Turing is interested in whether or not individual programs halt. I am interested in
trying to prove what are the bits, what is the numerical value, of the halting probability $\Omega$.
By the way, the value of $\Omega$
depends on your particular choice of programming language, which I don't
have time to discuss now.  $\Omega$ is also equal to the result of summing 1/2 
raised to powers which are the size in bits of every program that halts.
In other words, each $K$-bit program that halts contributes $1/2^K$
to $\Omega$.
    
And what precisely do I mean by an $N$-bit mathematical theory? Well, I'm thinking of
formal axiomatic theories, which are formulated using symbolic logic, not in any natural, human
language.  In such theories there are always a finite number of axioms and there are explicit
rules for mechanically deducing consequences of the axioms, which are called theorems.
An $N$-bit theory is one for which there is an $N$-bit program for systematically running
through the tree of all possible proofs deducing all the consequences of the axioms, which are
all the theorems in your formal theory. This is slow work, but in principle it can be done
mechanically, that's what counts.  David Hilbert believed that there had to be a single
formal axiomatic theory for all of mathematics; that's just another way of stating that
math is static and perfect and provides absolute truth.
    
Not only is this impossible, not only is Hilbert's dream impossible to achieve, but there
are in fact an infinity of irreducible mathematical truths, mathematical truths for which
essentially the only way to prove them is to add them as new axioms.  My first example of
such truths was determining elegant programs, and an even better example is provided by
the bits of $\Omega$. The bits of $\Omega$ are mathematical facts that are true for no reason
(no reason simpler than themselves), and thus violate Leibniz's principle of sufficient reason,
which states that if anything is true it has to be true for a reason.
    
In math the reason that something is true is called its proof. Why are the bits of $\Omega$
true for no reason, why can't you prove what their values are?
Because, as Leibniz himself points out in Sections 33 to 35 of \emph{The Monadology,}
the essence of the notion of proof is that you prove a complicated assertion by analyzing it,
by breaking it down until you reduce its truth to the truth of assertions that are so
simple that they no longer require any proof (self-evident axioms).
But if you cannot deduce the truth of something from any principle simpler than itself, then
proofs become useless, because \textbf{anything} can be proven from principles that are
equally complicated, e.g., by directly adding it as a new axiom without any proof.
And this is exactly what happens with the bits of $\Omega$.
    
\textbf{In other words, the normal, Hilbertian view of math is that all of mathematical truth,
an infinite number of truths, can be compressed into a finite number of axioms.
But there are an infinity of mathematical truths that cannot be compressed at all, not one bit!}
     
This is an amazing result, and I think that it has to have profound philosophical and
practical implications. Let me try to tell you why.
     
On the one hand, it suggests that pure math is more like biology than it is like physics.
In biology we deal with very complicated organisms and mechanisms, but in physics it is
normally assumed that there has to be a theory of everything, a simple set of equations
that would fit on a T-shirt and in principle explains the world, at least the physical world.
But we have seen that the world of mathematical ideas has infinite complexity, it cannot be
explained with any theory having a finite number of bits, which from a sufficiently abstract
point of view seems much more like biology, the domain of the complex, than like physics, 
where simple equations reign supreme.
    
On the other hand, this amazing result 
suggests that even though math and physics are different, they
may not be as different as most people think! I mean this in the following sense: 
In math you organize your computational experience, your lab is the computer, and in
physics you organize physical experience and have real labs.  But in both cases an
explanation has to be simpler than what it explains, and in both cases there are sets of
facts that cannot be explained, that are irreducible. Why? Well, in quantum physics it
is assumed that there are phenomena that when measured are equally likely to give either
of two answers (e.g., spin up, spin down) and that are inherently unpredictable and irreducible.
And in pure math we have a similar example, which is provided by the individual bits in 
the binary expansion of the numerical value of the halting probability $\Omega$.
     
This suggests to me a quasi-empirical view of math, in which one is more willing to 
add new axioms that are not at all self-evident but that are justified pragmatically, i.e.,
by their fruitful consequences, just like a physicist would.  I have taken the term
quasi-empirical from Lakatos. The collection of essays 
\emph{New Directions in the Philosophy of Mathematics}
edited by Tymoczko in my opinion pushes strongly in the direction of a quasi-empirical
view of math, and it contains an essay by Lakatos proposing the term ``quasi-empirical,'' as
well as essays of my own and by a number of other people.  Many of them may disagree
with me, and I'm sure do, but I repeat, in my opinion all of these essays justify
a quasi-empirical view of math, what I mean by quasi-empirical, which is somewhat
different from what Lakatos originally meant, but is in quite the same spirit, I think.
    
In a two-volume work full of important mathematical examples, 
Borwein, Bailey and Girgensohn have argued that experimental mathematics
is an extremely valuable research paradigm
that should be openly acknowledged and indeed vigorously embraced.
They do not go so far as to suggest that one should add new axioms
whenever they are helpful, without bothering with proofs, but they are certainly
going in that direction and nod approvingly at my attempts to provide some theoretical
justification for their entire enterprise by arguing that math and physics are not that different.
     
In fact, since I began to espouse these heretical views in the early 1970's, 
largely to deaf ears,
there have actually
been several examples of such new pragmatically justified, non-self-evident axioms: 
\begin{itemize}
\item
the $\mathbf{P \neq NP}$ hypothesis regarding
the time complexity of computations, 
\item
the axiom of projective determinacy in set theory,
and 
\item
increasing reliance on diverse unproved versions of the Riemann hypothesis regarding the distribution
of the primes.
\end{itemize}
So people don't need to have theoretical justification; they just do whatever is needed to get
the job done\ldots
    
The only problem with this computational and information-theoretic epistemology that I've 
just outlined to you is that it's based on the computer, 
and there are uncomputable reals.
So what do we do with contemporary physics which is full of partial differential equations and field
theories, all of which are formulated in terms of real numbers, \textbf{most of which are in fact
uncomputable}, as I'll now show.  
Well, it would be good to get rid of all that and convert to a \emph{digital physics.}
Might this in fact be possible?!  I'll discuss that too.

\section*{Computer Ontology: How real are real numbers?
What is the world made of?}

How did Turing prove that there are uncomputable reals in 1936? He did it like this.
Recall that the possible texts in French are a countable or denumerable infinity and can
be placed in an infinite list in which there is a first one, a second one, etc.
Now let's do the same thing with all the possible computer programs (first you have to
choose your programming language).  
So there is a first program, a second program, etc.
Every computable real can be calculated
digit by digit by some program in this list of all possible programs. 
Write the numerical value of that real next to the programs that calculate it,
and cross off the list all the programs that do not calculate an individual computable real.
We have converted a list of programs into a list of computable reals, and no computable
real is missing.
      
Next discard the integer parts of all these computable reals, and just keep the
decimal expansions.  Then put together a new real number by changing every digit on the
diagonal of this list (this is called Cantor's diagonal method; it comes from set
theory). So your new number's first digit differs from the first digit of the first
computable real, its second digit differs from the second digit of the second computable real,
its third digit differs from the third digit of the third computable real,
and so forth and so on.  
So it can't be in the list of all computable reals and it has to be uncomputable.
And that's Turing's uncomputable real number!\footnote
{\emph{Technical Note:} 
Because of \textbf{synonyms} like $.345999\ldots = .346000\ldots$
you should avoid having any 0 or 9 digits in Turing's number.}
    
Actually, there is a much easier way to see that there are uncomputable reals by using ideas
that go back to Emile Borel (again!). Technically, 
the argument that I'll now present uses what mathematicians call
\emph{measure theory,} which deals
with probabilities.  So let's just look at all the real numbers between 0 and 1.  These
correspond to points on a line, a line exactly one unit in length, whose leftmost point
is the number 0 and whose rightmost point is the number 1.  The total length of this
line segment is of course exactly one unit.  But I will now show you that all the
computable reals in this line segment can be covered using intervals whose total length
can be made as small as desired. In technical terms, the computable reals in the interval
from 0 to 1 are a set of measure zero, they have zero probability.
    
How do you cover all the computable reals? Well, remember that list of all the computable
reals that we just diagonalized over to get Turing's uncomputable real? This time let's
cover the first computable real with an interval of size $\epsilon/2$, let's cover
the second computable real with an interval of size $\epsilon/4$, and in general we'll cover
the $N$th computable real with an interval of size $\epsilon/2^N$.
The total length of all these intervals (which can conceivably overlap or fall partially 
outside  the unit interval from 0 to 1), is exactly equal to $\epsilon$, which can be
made as small as we wish!  In other words, there are arbitrarily small coverings, and
the computable reals are therefore a set of measure zero, 
they have zero probability, they constitute an infinitesimal
fraction of all the reals between 0 and 1.  So if you pick a real at random between 0 and 1,
with a uniform distribution of probability, it is infinitely unlikely, though possible, that
you will get a computable real!
    
What disturbing news! Uncomputable reals are not the exception, they are the majority!
How strange!
      
In fact, the situation is even worse than that. As Emile Borel points out 
on page 21 of his final book,
\emph{Les nombres inaccessibles} (1952), without making any reference to Turing, most
individual reals are not even uniquely specifiable, they cannot even be named or pointed out,
no matter how non-constructively, because of the limitations of human languages, which permit 
only a countable infinity of possible texts. The individually accessible  or nameable reals
are also a set of measure zero. Most reals are un-nameable, with probability one!
I rediscovered this result of Borel's on my own in a slightly different context, in which things can
be done a little more rigorously, which is when one is dealing with a formal axiomatic
theory or an \textbf{artificial} formal language instead of a natural human language.
That's how I present this idea in \emph{Meta Math!}.
    
So if most individual reals will forever escape us, why should we believe in them?!
Well, you will say, because they have a pretty structure and are a nice theory, a nice game
to play, with which I certainly agree, and also because they have important practical applications,
they are needed in physics. Well, perhaps not!  Perhaps physics can give up infinite precision
reals! How? Why should physicists want to do that?
    
Because it turns out that there are actually many reasons for being skeptical about the reals,
in classical physics, in quantum physics, and 
particularly
in more speculative contemporary efforts to
cobble together a theory of black holes and quantum gravity.
     
First of all, as my late colleague the physicist Rolf Landauer used to remind me, no
physical measurement has ever achieved more than a small number of digits of precision, not
more than, say, 15 or 20 digits at most, and such high-precision experiments are rare
masterpieces of the experimenter's art and not at all easy to achieve. 
    
This is only a practical limitation in classical physics. But in quantum physics
it is a consequence of the Heisenberg uncertainty principle and wave-particle duality
(de Broglie).  According to quantum theory, the more accurately you try to measure something,
the smaller the length scales you are trying to explore,
the higher the energy you need (the formula describing this involves Planck's constant).
That's why it is getting more and more expensive to build particle accelerators like
the one at CERN and at Fermilab, and governments are running out of money to fund high-energy
physics, leading to a paucity of new experimental data to inspire theoreticians.
     
Hopefully new physics will eventually emerge from astronomical observations of bizarre
new astrophysical phenomena, since we have run out of money here on earth! 
In fact, currently some of the most interesting physical speculations involve the thermodynamics 
of black holes, massive concentrations of matter that seem to be lurking at the hearts
of most galaxies. Work by Stephen Hawking and Jacob Bekenstein on the thermodynamics of
black holes suggests that any physical system can contain only a finite amount of
information, a finite number of bits whose possible maximum 
is determined by what is called the Bekenstein bound.
Strangely enough, this bound on the number of bits grows as the surface area of the physical
system, not as its volume, leading to the so-called ``holographic'' principle asserting
that in some sense space is actually two-dimensional even though it appears to have
three dimensions!
    
So perhaps continuity is an illusion, perhaps everything is really discrete.
There is another argument against the continuum
if you go down to what is called the Planck scale. At distances that extremely short our
current physics breaks down because spontaneous fluctuations in the quantum vacuum
should produce mini-black holes that completely tear spacetime apart. And that is not at all
what we see happening around us.
So perhaps distances that small \textbf{do not exist}.
    
Inspired by ideas like this, in addition to \emph{a priori} metaphysical biases in favor of
discreteness, a number of contemporary physicists have proposed building the world out of
discrete information, out of bits.  Some names that come to mind in this connection are
John Wheeler, Anton Zeilinger, Gerard 't Hooft, Lee Smolin, Seth Lloyd, Paola Zizzi, 
Jarmo M\"akel\"a and Ted Jacobson,
who are real physicists.
There is also more speculative work by a small cadre 
of cellular automata and computer enthusiasts
including Edward Fredkin and Stephen
Wolfram, whom I already mentioned, as well as Tommaso Toffoli, Norman Margolus, and others.
     
And there is also an increasing body of highly successful work on quantum computation and
quantum information that is not at all speculative, it is just a fundamental reworking
of standard 1920's quantum mechanics.   Whether or not quantum computers ever become practical,
the workers in this highly popular field have clearly established that it is illuminating
to study sub-atomic  quantum systems in terms of how they process qubits of 
quantum information and how
they perform computation with these qubits.
These notions have shed completely new light on the behavior of quantum mechanical
systems.
    
Furthermore, when dealing with complex systems such as those that occur in biology, 
thinking about information
processing is also crucial. As I believe Seth Lloyd said, the most important thing in 
understanding a complex system is to determine how it represents information and how
it processes that information, i.e., what kinds of computations are performed.
    
And how about the entire universe, can it be considered to be a computer? Yes,
it certainly can, it is constantly computing its future state from its current state,
it's constantly computing its own time-evolution!  And as I believe Tom Toffoli pointed out,
actual computers like your PC just hitch a ride on this universal computation!
    
So perhaps we are not doing violence to Nature by attempting to force her into a
digital, computational framework. Perhaps she has been flirting with us, giving us hints
all along, that she is really discrete, not continuous, hints that we choose not to hear,
because we are so much in love and don't want her to change!
    
For more on this kind of new physics, see the books by Smolin and von Baeyer in the bibliography.
Several more technical papers on this subject are also included there.

\section*{Conclusion}

Let me now wrap this up and try to give you a present to take home,
more precisely, a piece of homework.
In extremely abstract terms, I would say that the problem is, as was emphasized
by Ernst Mayr in his book \emph{This is Biology,} that the current philosophy of science deals
more with physics and mathematics than it does with biology.
But let me try to put this in more concrete terms and connect it with
the spine, with the central thread, of the ideas in this talk. 
    
To put it bluntly,
a closed, static, eternal fixed view of math can no longer be sustained.
As I try to illustrate with examples in my \emph{Meta Math!}\ book,
math actually advances by inventing new concepts, by completely
changing the viewpoint.
Here I emphasized new axioms, increased complexity, more information,
but what really counts are new ideas, new
concepts, new viewpoints.  And that leads me to the crucial question,
crucial for a proper open, dynamic, time-dependent view of mathematics,
\begin{center}
   \textbf{``Where do new mathematical ideas come from?''}
\end{center}
    
I repeat, math does not advance by mindlessly and mechanically 
grinding away deducing all the consequences of a fixed set of concepts
and axioms, not at all!   It advances with new concepts, new definitions,
new perspectives, through revolutionary change, paradigm shifts, not just by hard work.
    
In fact, I believe that this is actually the central question in biology
as well as in mathematics, it's the mystery of creation, of creativity: 
\begin{center}
   \textbf{``Where do new mathematical and biological ideas come from?''}
\\  
   \textbf{``How do they emerge?''}
\end{center}

Normally one equates a new biological idea with a new species, but in
fact every time a child is born, that's actually a new idea
incarnating; it's reinventing the notion of ``human being,'' which changes
constantly.
     
I have no idea how to answer this extremely important question; I wish I could.
Maybe \textbf{you} will be able to do it. Just try!
You might have to keep it cooking on a back burner while concentrating on other things, 
but don't give up!
All it takes is a new idea!
Somebody has to come up with it.
Why not you?\footnote
{I'm not denying the importance of Darwin's theory of evolution.
But I want much more than that, I want a profound, extremely general
mathematical theory that captures the essence of what life is and
why it evolves. I want a theory that gets to the heart of the matter.
And I suspect that any such theory will necessarily have to shed
new light on mathematical creativity as well.
Conversely, a deep theory of mathematical creation might also cover biological creativity.
    
A reaction from Gordana Dodig-Crnkovic:
``Regarding Darwin and Neo-Darwinism I agree with you --- it is a very good 
idea to go beyond. In my view there is nothing more beautiful and 
convincing than a good mathematical theory. And I do believe that it must 
be possible to express those thoughts in a much more general way\ldots\
I believe that it is a very 
crucial thing to try to formulate life in terms of computation. Not to say 
life is nothing more than a computation. But just to explore how far one 
can go with that idea. Computation seems to me a very powerful tool to 
illuminate many things about the material world and the material ground 
for mental phenomena (including creativity)\ldots\
Or would you suggest that creativity is given by God's will? That it is 
the very basic axiom?  
Isn't it possible to relate to pure chance? Chance and selection? Wouldn't 
it be a good idea to assume two principles: law and chance, where both are 
needed to reconstruct the universe in computational terms? (like chaos and 
cosmos?)''}

\section*{Appendix: Leibniz and the Law}

I am indebted to
Professor Ugo Pagallo for explaining to me that Leibniz, whose ideas and their elaboration
were the subject of my talk, is regarded as just as important in the field of
law as he is in the fields of mathematics and philosophy.
    
The
theme of my lecture was that if a law is arbitrarily complicated, then it is not a law;
this idea was traced via Hermann Weyl back to Leibniz.
In mathematics it leads to my $\Omega$ number and the surprising discovery of
completely lawless regions of mathematics, areas in which there is absolutely no
structure or pattern or way to understand what is happening.
    
The principle that an arbitrarily complicated law is not a law
can also be interpreted with reference to the legal system.
It is not a coincidence that the words
``law'' and ``proof'' and ``evidence'' are used in jurisprudence as well
as in science and mathematics.
In other words, the rule of law is equivalent to the rule of reason,
but if a law is sufficiently complicated,
then it can in fact be completely arbitrary and incomprehensible. 

\section*{Acknowledgements}

I wish to thank Gordana Dodig-Crnkovic for organizing E-CAP'05
and for inviting me to present the Turing lecture at E-CAP'05; also for stimulating
discussions reflected in those footnotes that are marked with GDC.
The remarks on biology are the product of a week spent in residence at
Rockefeller University in Manhattan, June 2005; I thank Albert Libchaber for
inviting me to give a series of lectures there to physicists and biologists.
The appendix is the result of lectures to philosophy of law students 
April 2005 at the Universities of Padua, Bologna and Turin;
I thank Ugo Pagallo for arranging this.
Thanks too to Paola Zizzi for help with the physics references.

\section*{References}

\begin{itemize}
\item
Edward Fredkin, http://www.digitalphilosophy.org.
\item
Stephen Wolfram, \emph{A New Kind of Science,} Wolfram Media, 2002.
\item
Gregory Chaitin, \emph{Meta Math!,} Pantheon, 2005.
\item
G. W. Leibniz, \emph{Discourse on Metaphysics, Principles of Nature and Grace, The Monadology,}
    1686, 1714, 1714.
\item
Hermann Weyl, \emph{The Open World,} Yale University Press, 1932.
\item
Thomas Tymoczko, \emph{New Directions in the Philosophy of Mathematics,} 
    Princeton University Press, 1998.
\item
Jonathan Borwein, David Bailey, Roland Girgensohn, 
    \emph{Mathematics by Experiment, Experimentation in Mathematics,} 
    A. K. Peters, 2003, 2004.
\item
Rebecca Goldstein, \emph{Incompleteness,} Norton, 2005. 
\item
B. Jack Copeland, \emph{The Essential Turing,} Oxford University Press, 2004. 
\item
Vladimir Tasic, \emph{Mathematics and the Roots of Postmodern Thought,} 
    Oxford University Press, 2001.
\item
Emile Borel, \emph{Les nombres inaccessibles,} Gauthier-Villars, 1952.
\item
Lee Smolin, \emph{Three Roads to Quantum Gravity,} Basic Books, 2001.
\item
Hans Christian von Baeyer, \emph{Information,} Harvard University Press, 2004.
\item
Ernst Mayr, \emph{This is Biology,} Harvard University Press, 1998.
\item
J. Wheeler, (the ``It from bit'' proposal), \emph{Sakharov Memorial Lectures on
Physics,} vol.\ 2, Nova Science, 1992.
\item
A. Zeilinger, (the principle of quantization of information), 
\emph{Found.\ Phys.}\ \textbf{29}:631--643 (1999).
\item
G. 't Hooft, ``The holographic principle,''
http://arxiv.org/hep-th/0003004.
\item
S. Lloyd, ``The computational universe,''
http://arxiv.org/quant-ph/0501135.
\item
P. Zizzi, ``A minimal model for quantum gravity,''
http://arxiv.org/gr-qc/0409069.
\item
J. M\"akel\"a, ``Accelerating observers, area and entropy,''
http://arxiv.org/gr-qc/0506087.
\item
T. Jacobson, ``Thermodynamics of spacetime,''
http://arxiv.org/gr-qc/9504004.
\end{itemize}

\end{document}